  %
  %


  \font \bbfive = bbm5
  \font \bbseven = bbm7
  \font \bbten = bbm10
  \font \eightbf = cmbx8
  \font \eighti = cmmi8 \skewchar \eighti = '177
  \font \eightit = cmti8
  \font \eightrm = cmr8
  \font \eightsl = cmsl8
  \font \eightsy = cmsy8 \skewchar \eightsy = '60
  \font \eighttt = cmtt8 \hyphenchar\eighttt = -1
  \font \sixbf = cmbx6
  \font \sixi = cmmi6 \skewchar \sixi = '177
  \font \sixrm = cmr6
  \font \sixsy = cmsy6 \skewchar \sixsy = '60
  \font \tensc = cmcsc10
  
  \font \titlefont = cmr7 scaled \magstep4
  \scriptfont \bffam = \bbseven
  \scriptscriptfont \bffam = \bbfive
  \textfont \bffam = \bbten


  \def \Headlines #1#2{\nopagenumbers
    \voffset = 2\baselineskip
    \advance \vsize by -\voffset
    \headline {\ifnum \pageno = 1 \hfil
    \else \ifodd \pageno \tensc \hfil \lowercase {#1} \hfil \folio
    \else \tensc \folio \hfil \lowercase {#2} \hfil
    \fi \fi }}

  \def \Title #1{\vbox{\baselineskip 20pt \titlefont \noindent #1}}

  \def \Date #1 {\footnote {}{\eightit Date: #1.}}

  \def \Authors #1{\bigskip \bigskip \noindent #1}

  \long \def \Addresses #1{\begingroup \eightpoint \parindent0pt
\medskip #1\par \par \endgroup }

  \long \def \Abstract #1{\begingroup \eightpoint
  \bigskip \bigskip \noindent
  {\sc ABSTRACT.} #1\par \par \endgroup }


  \def \goodbreak {\vskip0pt plus.1\vsize \penalty -250 \vskip0pt
plus-.1\vsize }

  \newcount \secno \secno = 0
  \newcount \stno

  \def \seqnumbering {\global \advance \stno by 1
    \number \secno .\number \stno }

  \def \label #1{\def\localsystemvariable {\number \secno
    \ifnum \number \stno = 0\else .\number \stno \fi }\global \edef
    #1{\localsystemvariable }}

  \def\section #1{\global\def\SectionName{#1}\stno = 0 \global
\advance \secno by 1 \bigskip \bigskip \goodbreak \noindent {\bf
\number \secno .\enspace #1.}\medskip \noindent \ignorespaces}

  \long \def \sysstate #1#2#3{\medbreak \noindent {\bf
\seqnumbering .\enspace #1.\enspace }{#2#3\vskip 0pt}\medbreak }
  \def \state #1 #2\par {\sysstate {#1}{\sl }{#2}}
  \def \definition #1\par {\sysstate {Definition}{\rm }{#1}}

  \def \proof {\medbreak \noindent {\it Proof.\enspace }}
  \def \proofend {\ifmmode \eqno \square \else \hfill \square
\looseness = -1 \medbreak \fi }

  \def \$#1{#1 $$$$ #1}

  \newcount \zitemno \zitemno = 0
  \def\numbertoi {{\rm(\ifcase \zitemno \or i\or ii\or iii\or iv\or
v\or vi\or vii\or viii\or ix\or x\or xi\or xii\fi)}}
  \def\zitemcntr {\global \advance \zitemno by 1 \numbertoi}
  \def\iItem #1{\medskip \item {#1}}
  \def\Item #1{\smallskip \item {#1}}
  \def\izitem {\zitemno = 0\iItem {\zitemcntr}}
  \def\zitem {\Item {\zitemcntr}}


  \def \({\left (}
  \def \){\right )}
  \def \[{\left \Vert }
  \def \]{\right \Vert }
  \def \*{\otimes }
  \def \+{\oplus }
  \def \:{\colon }
  \def \<{\left \langle }
  \def \>{\right \rangle }
  \def \text #1{\hbox {\rm #1}}
  \def \and {\hbox {,\quad and \quad }}
  
  \def \calcat #1{\,{\vrule height8pt depth4pt}_{\,#1}}

  \def \cstar {$C^*$}
  \def \for #1{,\quad #1}
  \def \inv {^{-1}}
  
  \def \square {\hbox {$\sqcap \!\!\!\!\sqcup $}}
  \def \stress #1{{\it #1}\/}

  \def \|{\Vert }
  \def \inv {^{-1}}

  \newskip \ttglue

  \def \eightpoint {\def \rm {\fam0 \eightrm }%
  \textfont0 = \eightrm
  \scriptfont0 = \sixrm \scriptscriptfont0 = \fiverm
  \textfont1 = \eighti
  \scriptfont1 = \sixi \scriptscriptfont1 = \fivei
  \textfont2 = \eightsy
  \scriptfont2 = \sixsy \scriptscriptfont2 = \fivesy
  \textfont3 = \tenex
  \scriptfont3 = \tenex \scriptscriptfont3 = \tenex
  \def \it {\fam \itfam \eightit }%
  \textfont \itfam = \eightit
  \def \sl {\fam \slfam \eightsl }%
  \textfont \slfam = \eightsl
  \def \bf {\fam \bffam \eightbf }%
  \textfont \bffam = \eightbf
  \scriptfont \bffam = \sixbf
  \scriptscriptfont \bffam = \fivebf
  \def \tt {\fam \ttfam \eighttt }%
  \textfont \ttfam = \eighttt
  \tt \ttglue = .5em plus.25em minus.15em
  \normalbaselineskip = 9pt
  \def \MF {{\manual opqr}\-{\manual stuq}}%
  \let \sc = \sixrm
  \let \big = \eightbig
  \setbox \strutbox = \hbox {\vrule height7pt depth2pt width0pt}%
  \normalbaselines \rm }


  \def\cite #1{{\rm [\bf #1\rm ]}}
  \def\scite #1#2{\cite{#1{\rm \hskip 0.7pt:\hskip 2pt #2}}}
  \def\lcite #1{#1}
  \def\bibitem#1#2#3#4{\smallskip \item {[#1]} #2, ``#3'', #4.}

  \def \references {
    \begingroup
    \bigskip \bigskip \goodbreak
    \eightpoint
    \centerline {\tensc References}
    \nobreak \medskip \frenchspacing }


  \def\L{{\cal L}}
  \def\B{{\cal B}}
  \def\H{{\cal H}}
  \def\LGH{L_2(G,\H)}
  \def\pil{\pi_\lambda}
  \def\d{\,d}
  \def\cf{{\it cf.}~}
  \def\sbe{_{B_e}}
  \def\bigarrow{{\hbox to 1cm{\rightarrowfill}}}
  
  \def\bx#1{\hbox to 150pt{\hfill $#1$ \hfill}}
  \def\comutriang#1{\vbox{
    \bx{\ C^*(\B)}
    \bx{^{\pil}\kern-3pt\swarrow \hfill \searrow^{\Lambda}}
    \bx{#1\hfill
      {\buildrel \Psi \over {\hbox to 1cm{\rightarrowfill}}}\hfill
      C^*_r(\B)}}}


  \def\Amena{E}
  \def\FD{FD}
  \def\JT{JT}
  \def\Ng{N}
  \def\Ped{P}
  \def\RieffelInduced{R}

  \Headlines
  {A Note on the Representation Theory of Fell Bundles} 
  {Ruy Exel} 

  \Date
  {April 1999}

  \Title
  {A Note on the Representation Theory of Fell Bundles}

  \Authors
  {Ruy Exel\footnote{*}{\eightrm Partially supported by CNPq.}}

  \Addresses
  {Departamento de Matem\'atica;
  Universidade Federal de Santa Catarina;
  88010-970 Florian\'opolis SC;
  Brazil.}

  \Abstract
  {We show that every Fell bundle $\B$ over a locally compact group
$G$ is \stress{proper} in a sense recently introduced by Ng.
Combining our results with those of Ng we show that if $\B$ satisfies
the \stress{approximation property} then it is amenable in the sense
that the full and reduced cross-sectional \cstar-algebras coincide.}

  \section{Introduction}
  Let $\B = \{B_t\}_{t\in G}$ be a Fell bundle over a locally compact
group $G$ (see \cite{\FD} for a comprehensive treatment of the theory
of Fell bundles, also referred to as \cstar-algebraic bundles).  We
denote by $L_2(\B)$ the right Hilbert $B_e$--module obtained by
completing the space $C_c(\B)$ of all continuous compactly supported
sections of $\B$, under the $B_e$--valued inner product
  (\cf \scite{\Amena}{Section 2}, \scite{\Ng}{2.2})
  given by
  $$
  \<\xi,\eta\>\sbe =
  \int \xi(t)^*\eta(t) \d t \for \xi,\eta\in C_c(\B).
  $$
  We should warn the reader that our notation for $L_2(\B)$ differs
from that used in \cite{\Ng}.

The \stress{left-regular representation} of $\B$ is the map (\cf
\scite{\Amena}{2.2})
  $$
  \Lambda : \B \to \L(L_2(\B)),
  $$
  where $\L(L_2(\B))$ indicates the \cstar-algebra of all adjointable
operators \scite{\JT}{1.1.7} on $L_2(\B)$, given for any $t$ in $G$
and any $b_t$ in $B_t$, by
  $$
  \Lambda(b_t) \xi \calcat s = b_t \xi(t\inv s)
  \for \xi\in L_2(\B), \quad s\in G.
  $$
  Let $C^*(\B)$ be the cross-sectional \cstar-algebra of $\B$ (\cf
\scite{\FD}{VIII.17.2}) defined to be the enveloping \cstar-algebra of
the Banach *-algebra $L_1(\B)$ formed by the integrable sections
\scite{\FD}{VIII.5.2}.
  The integrated form of $\Lambda$, which we also denote by $\Lambda$,
is the *-homomorphism
  $$
  \Lambda : C^*(\B) \to \L(L_2(\B))
  $$
  specified by setting $\Lambda(f) \xi = f * \xi$ (\cf
\scite{\Ng}{2.10}) for all $f$ in the dense subalgebra
$C_c(\B)\subseteq C^*(\B)$, and all $\xi\in C_c(\B)\subseteq L_2(\B)$.

  Suppose that we are given a *-representation (\cf
\scite{\FD}{VIII.8.2 and 9.1}) $\pi$ of $\B$ on a Hilbert space $\H$,
i.e, a map $\pi:\B\to\B(\H)$ that is linear on each fiber, that
satisfies
  \izitem $\pi(b)\pi(c) = \pi(bc)$,
  \zitem $\pi(b)^* = \pi(b^*)$,
  \medskip\noindent for each $b,c\in\B$, and that is
\stress{continuous} in the sense that for each $u\in\H$, the map
  $$
  b\in\B \mapsto \pi(b)u\in\H
  $$
  is continuous in the norm of $\H$.

  We may then form the representation $\pil$ of $\B$ on $L_2(G)\*\H =
\LGH$ by setting
  $$
  \pil(b_t) = \lambda_t\*\pi(b_t) \for t\in G\for b_t\in B_t,
  $$
  where $\lambda_t$ refers to the left-regular representation of $G$
on $L_2(G)$.
  We will also denote by
  $$
  \pil : C^*(\B) \to \B(\LGH)
  $$
  its integrated form \scite{\FD}{VIII.11.2, 11.4, 17.2}.

  Generalizing \scite{\Amena}{2.3}, Ng defines in \scite{\Ng}{2.11}
the \stress{reduced cross-sectional \cstar-algebra} of $\B$, denoted
$C^*_r(\B)$, to be $\Lambda(C^*(\B))$.  Ng also proposes an
alternative notion of reduced algebra, namely
  $$
  C^*_R(\B) := \pil(C^*(\B)),
  $$
  where $\pi$ is any faithful *-representation of $\B$ on a Hilbert
space $\H$.

  There exists (see below) a unique surjective *-homomorphism
  $\Psi: C^*_R(\B) \to C^*_r(\B)$ such that the diagram
  $$
  \comutriang{C^*_{R}(\B)}
  $$
  commutes.  Ng thus introduced the notion of \stress{proper} Fell
bundles (\cf \scite{\Ng}{2.15}) to single out those for which $\Psi$ is
injective.  It is noticed in \cite{\Ng} that Theorem 3.3 in
\cite{\Amena} implies that Fell bundles over discrete groups are
automatically proper.  It is also shown that saturated Fell bundles
are always proper \scite{\Ng}{2.17}, as well as those whose underlying
group is compact \scite{\Ng}{A.3}.

  It is the purpose of this note to show that all Fell bundles are
proper and hence that the alternative reduced algebra $C^*_R(\B)$
proposed by Ng always coincides with $C^*_r(\B)$.

One of the main consequences is that the properness hypothesis
required in the main result of \cite{\Ng} (Proposition 3.9) becomes
superfluous and hence we conclude that all Fell bundles satisfying the
\stress{approximation property} (Definition 3.6 in \cite{\Ng}; see
also \scite{\Amena}{4.4}) are amenable in the sense that $\Lambda$ is
an isomorphism from $C^*(\B)$ to $C^*_r(\B)$.

  \section{Preliminaries}
  Let us fix, throughout, a *-representation $\pi:\B\to\B(\H)$.
  Restricting $\pi$ to $B_e$ we may view $\H$ as a left $B_e$--module
and hence we may form the tensor product
  $L_2(\B)\*\sbe\H$
  (\cf \scite{\RieffelInduced}{5.1}, \scite{\JT}{1.2.3}), which is a
Hilbert space under the inner product defined by
  $$
  \<\xi\*u,\eta\*v\> =
  \<u,\pi\(\<\xi,\eta\>\sbe\) v\>
  \for \xi,\eta\in L_2(\B)\for u,v\in\H.
  $$

  \state Proposition
  \label \PropoV
  {\rm(Lemma 2.4 in \cite{\Ng})}.
  There exists an isometry
  $$
  V: L_2(\B)\*\sbe\H \to \LGH,
  $$
  such that for all $\xi\in L_2(\B)$, $u\in\H$, and $t\in G$ one has
  $$
  V(\xi\*u)\calcat t = \pi(\xi(t))u.
  $$

  \proof
  It is obvious that $V$ is balanced with respect to the corresponding
actions of $B_e$ and hence it is well defined on the algebraic tensor
product
  $L_2(\B)\odot\sbe\H$.  Now let
  $\xi,\eta\in L_2(\B)$, and $u,v\in\H$.  We have
  $$
  \<V(\xi\*u),V(\eta\*v)\> =
  \int \<\pi(\xi(t))u,\pi(\eta(t))v\> \d t =
  \<u,\pi\(\int\xi(t)^*\eta(t)\d t\)v\> \$=
  \<u,\pi(\<\xi,\eta\>\sbe)v\> =
  \<\xi\*u,\eta\*v\>,
  $$
  from which all of the remaining details follow.
  \proofend

  It should be observed that $V$ is not necessarily surjective.  In
fact, note that the vector $V(\xi\*u)\calcat t$, mentioned above, lies
in $\pi(B_t)\H$ which is often a proper subset of $\H$.  This is
related to the notion of \stress{saturated} representations
\scite{\Ng}{Definition 2.5} and is one of the main stumbling blocks we
must overcome in order to achieve our goals.

  \state Proposition
  \label \TensorId
  {\rm(Lemma 1.3 in \cite{\Ng})}.
  If $\pi|\sbe$ is injective then so is
  the *-homomorphism
  $$
  T \in \L(L_2(\B)) \longmapsto T\*1 \in \B(L_2(\B)\*\sbe\H).
  $$

  \proof Suppose that $T\*1=0$.  Then, for all
  $\xi,\eta\in L_2(\B)$, and $u,v\in\H$ we have
  $$
  0 =
  \<(T\*1)(\xi\*u),\eta\*v\> =
  \<u,\pi\(\<T(\xi),\eta\>\sbe\) v\>.
  $$
  Since $u$ and $v$ are arbitrary, and $\pi$ is supposed injective on
$B_e$, this implies that $\<T(\xi),\eta\>\sbe=0$ for all $\xi$ and
$\eta$, which in turn gives $T=0$.
  \proofend

  In particular, when $\pi|\sbe$ is injective, we have by
\lcite{\TensorId} that $C^*_r(\B)$ is isomorphic to the algebra
$\Lambda(C^*(\B))\*1$ of operators on the Hilbert space
$L_2(\B)\*\sbe\H$.

  \state Proposition
  \label \SquareDiagram
  For any $b\in \B$
  the diagram
  $$
  \matrix{
  &{\scriptstyle \Lambda(b)\*1} \cr
  L_2(\B)\*\sbe\H & \bigarrow & L_2(\B)\*\sbe\H \cr\cr
  {\scriptstyle V}\Big\downarrow&& \Big\downarrow {\scriptstyle
V}\cr\cr
  \LGH & \bigarrow & \LGH \cr
  &{\scriptstyle \pil(b)}\cr
  }
  $$
  commutes.

  \proof
  Let $t\in G$ be such that $b\in B_t$.  We then have for all $\xi\in
L_2(B)$, $u\in\H$, and $s\in G$ that
  $$
  V(\Lambda(b)\*1)(\xi\*u)\calcat s =
  V\big(\Lambda(b)\xi\*u\big)\calcat s =
  \pi\(\Lambda(b)\xi\calcat s\) u =
  \pi(b\xi(t\inv s))u.
  $$
  On the other hand
  $$
  \pil(b)V(\xi\*u)\calcat s =
  \pi(b)\(V(\xi\*u)\calcat {t\inv s}\) =
  \pi(b)\pi(\xi(t\inv s))u.
  \proofend
  $$

  It follows that the same holds if, in place of the ``$b$'' in the
statement above, we substitute any $a\in C^*(\B)$, since the
corresponding representations at the level of $C^*(\B)$ are integrated
from those of $\B$.

  \definition
  (\cf \cite{\Ng}).
  Given a *-representation $\pi: \B \to \H$ as above we shall denote
by $C^*_{R,\pi}(\B)$ the algebra $\pil(C^*(\B))$ of operators on
$\LGH$.

  When $\pi|\sbe$ is faithful, $C^*_{R,\pi}(\B)$ was proposed by Ng
\cite{\Ng} as an alternative reduced cross-sectional \cstar-algebra
for $\B$.
  The first relationship between $C^*_{R,\pi}(\B)$ and $C^*_r(\B)$ is
given by:

  \state Proposition
  \label \Triangle
  Suppose that $\pi|\sbe$ is injective.  Then
  for any $a\in C^*(\B)$ one has that
  $
  \[\Lambda(a)\] \leq
  \[\pil(a)\].
  $
  Therefore there exists a unique *-homomorphism
  $\Psi: C^*_{R,\pi}(\B) \to C^*_r(\B)$ such that the diagram
  $$
  \comutriang{C^*_{R,\pi}(\B)}
  $$
  commutes.

  \proof
  By \lcite{\SquareDiagram} we have that $\Lambda\*1$ is equivalent to
a subrepresentation of $\pil$.  Therefore
  $$
  \[\Lambda(a)\* 1\] \leq
  \[\pil(a)\].
  $$
  Now, by \lcite{\TensorId}, we have that
  $ \[\Lambda(a)\* 1\] = \[\Lambda(a)\]$.
  The existence of $\Psi$ now follows by routine arguments.
  \proofend

  \section{The main result}
  As already indicated, we plan to prove that $\Psi$ is an isomorphism
under the hypothesis that $\pi|\sbe$ is injective.  This is clearly
equivalent to proving that for any $a\in C^*(\B)$ one has that
  $
  \[\Lambda(a)\] =
  \[\pil(a)\].
  $
  The starting point is that, although $\Lambda\*1$ is but a
subrepresentation of $\pil$, we may ``move it around'' filling out the
whole of the representation space for $\pil$.  What will do the
``moving around'' will be the \stress{right-regular representation} of
$G$, namely the unitary representation $\rho$ of $G$ on $L_2(G)$ given
by
  $$
  \rho_r(\xi)\calcat s = \Delta(r)^{1/2}\xi(sr),
  $$
  for $\xi\in L_2(G)$, and $r,s\in G$, where $\Delta$ is, as usual,
the modular function for $G$.

  \state Proposition
  \label \RhoCommute
  For each $r\in G$,
  \izitem The unitary operator $\rho_r\*1$ on $L_2(G)\*\H = \LGH$
lies in the commutant of $\pil(C^*(\B))$.
  \zitem Consider the isometry
  $$
  V_r: L_2(\B)\*\sbe\H \to \LGH,
  $$
  given by $V_r = (\rho_r\*1) V$.  Then for all $a\in C^*(\B)$ one
has $V_r (\Lambda(a)\*1) = \pil(a) V_r.$
  \zitem Let $K_r$ be the range of\/ $V_r$.  Then $K_r$ is invariant
under $\pil$ and the restriction of\/ $\pil$ to $K_r$ is equivalent to
$\Lambda\*1$.

  \proof 
  It is clear that $\rho_r\*1$ commutes with $\pil(b_t) =
\lambda_t\*\pi(b_t)$ for any $b_t\in B_t$.  It then follows that
$\rho_r\*1$ also commutes with the range of the integrated form of
$\pil$, whence (i).  The second point follows immediately from (i) and
\lcite{\SquareDiagram}. Finally, (iii) follows from (ii).
  \proofend

  Our next result is intended to show that the $K_r$'s do indeed fill
out the whole of $\LGH$.

  \state Proposition
  \label \Density
  Suppose that $\pi|\sbe$ is non-degenerate.  Then
  the linear span of\/ $\bigcup_{r\in G} K_r$ is dense in $\LGH$.

  \proof Let
  $$
  \Gamma = {\rm span}\{V_r(\xi\*u): r\in G,\ \xi\in C_c(\B),\ u\in
\H\}.
  $$
  Since
  $$
  V_r(\xi\*u)\calcat t =
  (\rho_r\*1)V(\xi\*u)\calcat t =
  \Delta(r)^{1/2}V(\xi\*u)\calcat {tr} =
  \Delta(r)^{1/2}\pi(\xi(tr))u
  \for t\in G,
  $$
  and since we are taking $\xi$ in $C_c(\B)$ above, it is easy to see
that $\Gamma$ is a subset of $C_c(G,\H)$.  Our strategy will be to use
  \scite{\FD}{II.15.10} for which we must prove that:
  \medskip
  \item{(I)} If $f$ is a continuous complex function on $G$ and
$\eta\in\Gamma$, then the pointwise product $f\eta$ is in $\Gamma$;
  \item{(II)} For each $t\in G$ the set $\{\eta(t):\eta\in \Gamma\}$
is dense in $\H$.

  \medskip
  The proof of (I) is elementary in view of the fact that $C_c(\B)$ is
closed under pointwise multiplication by continuous scalar-valued
functions \scite{\FD}{II.13.14}.
  In order to prove (II) let $v\in\H$ have the form
  $v=\pi(b)u$, where $b\in B_e$ and $u\in\H$.  By
\scite{\FD}{II.13.19} let $\xi\in C_c(\B)$ be such that $\xi(e)=b$.
It follows that
  $\eta_r:=V_r(\xi\*u)$ is in $\Gamma$ for all $r$.
  Also note that ,
  setting $r=t\inv$, we have
  $$
  \eta_{t\inv}(t) =
  \Delta(t)^{-1/2}\pi(\xi(e))u =
  \Delta(t)^{-1/2}\pi(b)u =
  \Delta(t)^{-1/2}v.
  $$
  This shows that $v\in\{\eta(t):\eta\in \Gamma\}$.  Since the set of
such $v$'s is dense in $\H$, because $\pi|\sbe$ is non-degenerate,
we have that (II) is proven.

  As already indicated, it now follows from \scite{\FD}{II.15.10} that
$\Gamma$ is dense in $\LGH$.  Since $\Gamma$ is contained in the
linear span of $\bigcup_{r\in G} K_r$, the conclusion follows.
  \proofend

  The following is our main technical result:

  \state Lemma
  \label \Mainlemma
  For all $a\in C^*(\B)$ one has that $\[\pil(a)\] \leq
\[\Lambda(a)\]$.

  \proof We may clearly suppose, without loss of generality, that
$\pi$ is non-degenerate.  By \scite{\FD}{VIII.9.4} it follows that
$\pi|\sbe$ is non-degenerate as well.
  Under this assumption we claim that for all $a\in C^*(\B)$ one has
that
  $$
  \Lambda(a) = 0 \quad =\!\!\Rightarrow \quad \pil(a) = 0.
  $$
  In order to see this suppose that $\Lambda(a) = 0$.
  Then for each $r\in G$ we have by \lcite{\RhoCommute}.(ii) that
  $
  \pil(a) V_r = V_r (\Lambda(a)\*1) =0.
  $
  Therefore
  $\pil(a)=0$ in the range $K_r$ of $V_r$.
  By \lcite{\Density} it folows that $\pil(a)=0$, thus proving our
claim.

  Define a map
  $$
  \varphi : C^*_r(\B) \longrightarrow \B(\LGH)
  $$
  by
  $\varphi (\Lambda(a)) := \pil(a)$, for all $a$ in $C^*(\B)$.  By the
claim above we have that $\varphi$ is well defined.  Also, it is easy
to see that $\varphi$ is  a *-homomorphism.  It follows that
$\varphi$ is contractive and hence that for all $a$ in $C^*(\B)$
  $$
  \[\pil(a)\] =
  \[\varphi (\Lambda(a))\] \leq
  \[\Lambda(a)\].
  \proofend
  $$

  Our main results follow more or less immediately from
\lcite{\Mainlemma}:

  \state Corollary
  \label \AlwaysProper
  Let $\B$ be any Fell bundle over a locally compact group $G$ and let
$\pi$ be a *-representation of $\B$ on the Hilbert space $\H$ such
that $\pi|\sbe$ is injective.  Then the map
  $\Psi: C^*_{R,\pi}(\B) \to C^*_r(\B)$
  defined above is an isomorphism.  Therefore $\B$ is always proper in
the sense of Ng \cite{\Ng}.

  \state Corollary
  Suppose that the Fell bundle $\B$ satisfies the approximation
property (Definition 3.6 in \cite{\Ng}; see also \scite{\Amena}{4.4}),
then $\B$ is amenable in the sense that $\Lambda$ is an isomorphism
from $C^*(\B)$ to $C^*_r(\B)$.

  \proof Combine Proposition 3.9 in \cite{\Ng} with
\lcite{\AlwaysProper}.
  \proofend

  The following generalizes \scite{\Ped}{7.7.5} to the context of Fell
bundles:

  \state Corollary
  Let $\pi:\B\to\B(\H)$ be a representation of the Fell bundle $\B$
and let $\pil$ be the representation of $\B$ on $\LGH$ given by
  $\pil(b_t) = \lambda_t\*\pi(b_t)$, for $t\in G$, and $b_t\in B_t$.
Denote also by $\pil$ the representation of $C^*(\B)$ obtained by
integrating $\pil$.  Then $\pil$ factors through $C^*_r(\B)$.
Moreover, in case $\pi|\sbe$ is faithful, the representation of
$C^*_r(\B)$ arising from this factorization is also faithful.

  \proof
  Follows immediately from \lcite{\Triangle} and \lcite{\Mainlemma}.
  \proofend

  \endgroup

  \bye